\documentclass{amsart}

\usepackage{amsmath}%
\usepackage{amsfonts}%
\usepackage{amssymb}%
\usepackage{graphicx}
\newtheorem{theorem}{Theorem}
\theoremstyle{plain}

\newtheorem{corollary}{Corollary}

\newtheorem{definition}{Definition}

\newtheorem{lemma}{Lemma}

\newtheorem{proposition}{Proposition}
\newtheorem{remark}{Remark}

\numberwithin{equation}{section}
\newcommand{\R}{\mathbb{R}}
\newcommand{\C}{\mathbb{C}}
\DeclareMathOperator{\imag}{Im}
\DeclareMathOperator{\ind}{Ind}

\DeclareMathOperator{\spn}{span}
\DeclareMathOperator{\Lie}{Lie}
\begin{document}
\title[Global representations of the Heat and Schr\"{o}dinger equation]{Global representations of the Heat and Schr\"{o}dinger equation with singular potential}

\author{Jose A. Franco}  
\address{University of North Florida \\ 1 UNF Drive, Jacksonville, FL 32082 }
\email{jose.franco@unf.edu}

\author{Mark R. Sepanski}
\address{Baylor University \\ One Bear Place \# 97328, Waco, TX 76798  }
\email{mark\_sepanski@baylor.edu}


\begin{abstract}
We study the $n$-dimensional Schr\"{o}dinger equation with a singular potential $V_\lambda(x)=\lambda \left\|x\right\|^{-2}$. Its solution space is studied as a global representation of $\widetilde{SL(2,\R)}\times O(n)$. A special subspace of solutions for which the action globalizes is constructed via nonstandard induction outside the semisimple category. The space of $K$-finite vectors is calculated, obtaining conditions for $\lambda$ so that this space is non-empty. The direct sum of solution spaces over such admissible values of $\lambda$ is studied as a representation of the $2n+1$-dimensional Heisenberg group.
\end{abstract}

\maketitle

\section{Introduction}

The Schr\"{o}dinger and heat equations have been heavily studied in physics and in mathematics. In physics, the study of the Schr\"{o}dinger equation with  inverse square potential is important in the study of the motion of a dipole in a cosmic string background, as noted by Bouaziz and Bawin \cite{Bouaziz}. Galajinsky, Lechtenfeld, and Polovnikov\cite{Galajinsky} studied it in the context of the Calogero model of a set of decoupled particles on the real line. This potential is also relevant in the fabrication of nanoscale atom optical devices, the study of dipole-bound anions of polar molecules, and in the study of the behavior of three-body systems in nuclear physics (see Bawin and Coon \cite{Bawin}). A generalization of this potential was used by Cambalong, Epele, Fanchiotti, and Canal \cite{Horacio} to study the symmetries of the interaction of a polar molecule and an electron. The inverse square potential is also relevant in the study of the Efimov spectrum of the three body problem\cite{Efimov} among 
other applications. In this paper, we calculate a basis for the space of solutions to the Schr\"{o}dinger equation on which the action of the symmetry group is completely described.

In mathematics, these equations have been studied from many points of view. One of these points of view is the symmetry group analysis associated to them. Along this line,  Sepanski and Stanke \cite{Sepanski1} examined the solutions to a family of differential equations that includes the potential free Schr\"{o}dinger and heat equations, as global representations of the corresponding Lie symmetry group, $G:=(\widetilde{SL(2,\R)}\times O(n) )\ltimes H_{2n+1}$, where $\widetilde{SL(2,\R)}$ denotes the two-fold cover of $SL(2,\R)$, $O(n)$ denotes  the orthogonal group, and $H_{2n+1}$ denotes the $2n+1$-dimensional Heisenberg group. Using the same techinques, one of the authors \cite{Franco} studied an invariant subspace of solutions to the one-dimensional Schr\"{o}dinger equation with singular potential $V(x)=\lambda/x^2$. The natural generalization of these two articles is the goal of this article. 

Let $x \in \R^n$, $s \in \C$, and $\Delta_n$ denote the Laplacian operator on $\R^n$.  Here we study the representation theory associated to special subspaces of solutions of the family of differential equations 
\begin{equation}\label{MainProblem}4s\partial_t+\Delta_n=\frac{2\lambda}{\left\|x\right\|^2}\end{equation}
that are invariant under the action of the group $\widetilde{SL(2,\R)}\times O(n)$ or $G$ when $\lambda=0$.  By letting $s=i/2$ or $s=-1/4$ one obtains the Schr\"{o}dinger or the heat equation with singular potential, respectively. 

In terms of representation theory, this problem will be equivalent to the solution of an eigenvalue problem associated to a Casimir element in a certain line bundle over a compactification of $\R^{1,n}$. In more detail, we start by constructing an induced representation space that carries the structure of a global Lie group representation. To do that, we consider a parabolic-like subgroup $\overline{P}$ of $G$ and the smoothly induced representation $$I(q,r,s)=\ind_{\overline P}^{G}(\chi_{q,r,s})$$ where $\chi_{q,r,s}$ is a character on $\overline{P}$ with parameters $r,s \in \C$ and $q \in \mathbb Z_4$ (See Section \ref{IndRepsSection}). Using the fact that $\R^{1,n}$ can be embedded as an open dense set in $G/\overline{P}$, restriction  to  $G/\overline{P}$ can be  used to realize $I(q,r,s)$ as a space of smooth functions with certain decay conditions. We denote this space by $I'(q,r,s)$. In the semisimple category, this is called the non-compact picture.

Let $\Omega'$ denote the action of $2\Omega_{\mathfrak{sl}_(2,\R)}-\Omega_{\mathfrak{so}(n)}-r(r+2)$ on $I'(q,r,s)$, where $\Omega_{\mathfrak{sl}(2,\R)}$ is the Casimir element of $ \mathfrak{sl}(2,\R)$ and $\Omega_{\mathfrak{so}(n)}$ is the Casimir element for $\mathfrak{so}(n)$ in the universal enveloping algebra of $\Lie(G)$.  We show that the kernel of $\Omega'-2\lambda$ in $I'(q,r,s)$ is a group invariant subspace of solutions to \eqref{MainProblem} on $\R^{1,n}$.

In order to study the structure of $\ker(\Omega'-2\lambda)$, we turn to the analog of what would be the compact picture in the semisimple category. There we explicitly derive the conditions on $\lambda$ for the existence of non-trivial $K$-finite vectors (see Theorem \ref{KTypesForm}). When they exist, the form of the $K$-finite vectors is given explicitly in terms of confluent hypergeometric functions of the first kind and harmonic polynomials. 

We will say that the eigenvalue $\lambda$ is admissible if and only if the space of $K$-finite vectors, $\ker(\Omega'-2\lambda)_K$, is non-trivial. The set of all non-zero admissible eigenvalues will be denoted by $A_n$ and is explicitly determined in Section \ref{Struc2}. For each $\lambda \in A_n$, $\ker(\Omega'-2\lambda)_K$ decomposes under $\mathfrak{sl}_2\times O(n)$ as a direct sum of finitely many infinite dimensional representations. The structure of the modules is completely determined including the determination of when highest or lowest weights exist.

When $\lambda =0$, $\ker \Omega'$ is also invariant under the action of the Heisenberg algebra and its structure is determined in \cite{Sepanski1}. However, for non-zero values of $\lambda$, the action of the Heisenberg algebra does not preserve $\ker(\Omega'-2\lambda)_K$. Nevertheless, we show that the space $$\ker\Omega' \oplus \bigoplus_{\lambda\in A_n}\ker(\Omega'-2\lambda)_K$$ carries the structure of a $\mathfrak{g}$-module, where $\mathfrak{g}$ denotes the Lie algebra of $G$. The composition series of this space is determined in Theorem \ref{BigTheorem}.

\section{Notation}
We will follow the constructions in Sepanski and Stanke\cite{Sepanski1} very closely.
\subsection{The Group}

For $x, y \in \R^{2n}$, let $\left\langle x,y\right\rangle=x^TJ_ny$ where $J_n=\left(\begin{smallmatrix}0 & I_n  \\ -I_n & 0\end{smallmatrix}\right)$. Let $H_{2n+1}$ denote the $(2n+1)$-dimensional Heisenberg group. The multiplication in $H_{2n+1}$ is given by $$(v,t)(v',t')=(v+v',t+t'+\langle v,v'\rangle)$$ where $v,v' \in \R^{2n}$ and $t,t' \in \R$. 

An element $\sigma \in Sp(2n,\R)$ acts on $H_{2n+1}$ by $\sigma.(v,t)=(\sigma.v,t)$ where the action $\sigma.v$ is the standard action of $Sp(2n,\R)$ on $\R^{2n}$. Thus we can define the product on $Sp(2n,\R)\ltimes H_{2n+1}$ by $$(\sigma, h)(\tau,k)=(\sigma \tau, \tau^{-1}(h)k)$$ for  $\sigma, \tau \in Sp(2n,\R)$ and $h,k \in H_{2n+1}$.

The group $SL(2,\R)=Sp(2,\R)$ can be embedded in $Sp(2n,\R)$ by $\left(\begin{smallmatrix} a & b \\ c & d \end{smallmatrix}\right)\mapsto \left(\begin{smallmatrix} aI_n & bI_n \\ cI_n & dI_n \end{smallmatrix}\right)$ and the group $O(n)$ can be embedded diagonally by $u \mapsto \left(\begin{smallmatrix} u & 0 \\ 0 & u \end{smallmatrix}\right)$. Since these two images commute, there exists a homomorphism $B:SL(2,\R)\times O(n) \to Sp(2n,\R)$ with kernel $\pm (I_2\times I_n)$.

Following the realization of the two-fold cover of $SL(2,\R)$ of Kashiwara and Vergne \cite{Kashi}, define the complex upper half plane $D:=\{z\in\C|\imag z >0\}$ and let $SL(2,\R)$ act on $D$ by linear fractional  transformations, that is, if $g=\left(\begin{smallmatrix} a & b \\ c & d \end{smallmatrix}\right)\in SL(2,\R)$ and $z\in D$ then $$g.z=\frac{az+b}{cz+d}.$$ Define $d:SL(2,\R)\times D \to \C$ by $d(g,z):=cz+d$. Then there are exactly two smooth square roots of $d(g,z)$ for each $g \in SL(2,\R)$ and $z\in D$. The double cover can be realized as:
\begin{multline*}
G_2:=\left\{(g,\epsilon) | g\in SL(2,\R) \text{ and smooth }\epsilon:D\to\C \right. \\ \left. \text{ such that } \epsilon(z)^2=d(g,z)\text{ for } z\in D \right\}
\end{multline*}
with the product defined by $$(g_1,\epsilon_1(z))(g_2,\epsilon_2(z))=(g_1g_2, \epsilon_1(g_2.z)\epsilon_2(z)).$$

Let $p:G_2\to SL(2,\R)$ be the canonical projection. Then $B\circ (p\otimes 1):G_2\times O(n) \to Sp(2n,\R)$ is a homomorphism and the semidirect product $$G := (G_2\times O(n) )\ltimes H_{2n+1}$$ is well-defined via this homomorphism.

\subsection{Parabolic Subgroup and Induced Representations}\label{IndRepsSection}

Let  $\exp_{G_2}: \mathfrak{sl}(2,\R)\to G_2$ denote the exponential map. Let $\mathfrak a=\{\left(\begin{smallmatrix}t & 0 \\ 0&-t\end{smallmatrix}\right)| t\in\R\}$, $\mathfrak n=\{\left(\begin{smallmatrix}0 & t \\ 0&0\end{smallmatrix}\right)| t\in\R\}$, and $\overline{\mathfrak n}=\{\left(\begin{smallmatrix}0 & 0 \\ t&0\end{smallmatrix}\right)| t\in\R\}$. Then, 
\begin{gather*}
A:= \exp_{G_2}(\mathfrak a)=\{(\left(\begin{smallmatrix}t & 0 \\ 0&t^{-1}\end{smallmatrix}\right),z\mapsto e^{-t/2})| t\in\R^{\geq 0}\}\\
N:= \exp_{G_2}(\mathfrak n)=\{(\left(\begin{smallmatrix}1 & t \\ 0& 1\end{smallmatrix}\right), z\mapsto 1)|t\in\R\}\\
\overline{N}:= \exp_{G_2}(\overline{\mathfrak{n}})=\{(\left(\begin{smallmatrix}1 & 0 \\ t & 1\end{smallmatrix}\right), z\mapsto \sqrt{tz+1})|t\in\R\}.
\end{gather*}
Let $\mathfrak{k}:= \{\left(\begin{smallmatrix}0 & \theta \\ -\theta & 0 \end{smallmatrix}\right):\theta \in \R\}$ then $$K_2:= \exp_{G_2}(\mathfrak k)=\{(g_\theta,\epsilon_\theta):= (\left(\begin{smallmatrix}\cos\theta & \sin\theta \\ -\sin\theta & \cos\theta\end{smallmatrix}\right),z \mapsto \sqrt{\cos\theta-z\sin\theta} )|\theta \in \R\}$$
where $\sqrt{\cdot}$ denotes the principal square root in $\C$.  Writing $M$ for the centralizer of $A$ in $K_2$ then
$$M=\{m_j:=(\left(\begin{smallmatrix}-1 & 0 \\0 &-1\end{smallmatrix}\right)^j, z\to i^{-j}) | j=0,1,2,3\}.$$
Let $W \subset H_{2n+1}$ be given by $W=\{(0,v,w)| v\in \R^n \text{ and }w \in \R \}\cong \R^{n+1}$ and let $\overline P=(MA\overline N \times O(n))\ltimes W$.

It is well known that the character group on $A$ is isomorphic to the additive group $\C$ so any character on $A$ can be indexed by a constant $r \in \C$ and defined by 
$$\chi_r\bigl((\left(\begin{smallmatrix}t & 0 \\ 0&t^{-1}\end{smallmatrix}\right),z\mapsto e^{-t/2})\bigr)=t^r$$
for $t>0$. A character on $M$ is parametrized by $q\in \mathbb Z_4$ and defined by $\chi_q(m_j)=i^{j q}$. A character on $W$ can be parametrized by $s\in \C$ and defined by,
$$\chi_{s}\bigl((0,v,w)\bigr)=e^{sw}.$$
 Finally, any character on $\overline P$ that is trivial on $N$ is parametrized by a triplet $(q,r,s)$ where $s, r \in \C$ and $q\in \mathbb Z_4$ and defined by
\begin{equation}\label{Character}\chi_{q,r,s}\left(((-1)^j\left(\begin{smallmatrix} a& 0 \\ c & a^{-1}\end{smallmatrix}\right),z\mapsto i^{-j}e^{-a/2}\sqrt{acz+1}),\left( 
0,v,w\right) \right)=i^{jq}|a|^r e^{sw}.\end{equation}
The representation space induced by $\chi_{q,r,s}$ will be denoted by $I(q,r,s)$ and defined by 
$$I(q,r,s):= \{\phi:G\to\C | \phi \in C^\infty \text{ and } \phi(g\overline p)=\chi_{q,r,s}^{-1}(\overline p)\phi(g) \text{ for } g\in G, \overline p\in \overline P\}.$$
The action of $G$ on $I(q,r,s)$ is given by left translation: $(g_1.\phi)(g_2)=\phi(g_1^{-1}g_2)$.

\section{The Non-compact Picture}

If $X:=\{(x,0,0) |x \in \R^n\}$, then $H_{2n+1}=XW$ and $(N\times X)\overline P$ is open and dense in $G$. Since $N\times X$ is isomorphic to $\R^{n+1}$ via $(t,x)\mapsto N_{t,x}:=\left[(\left(\begin{smallmatrix}1 & t \\ 0& 1\end{smallmatrix}\right), z\mapsto 1),(x,0,0)\right]$ and since a section in the induced representation, $I(q,r,s)$, is determined by its restriction to $N$, there exists an injection of $I(q,r,s)$ into $C^\infty(\R^{n+1})$, given by restriction of domain. The image of this map is identified as
\begin{equation*}I'(q,r,s)=\{f\in C^\infty(\R^{n+1}) | f(t,x)=\phi(N_{t,x})  \text{ for some } \phi\in I(q,r,s)\}\end{equation*} 
and is given the $G$-module that makes the map $\phi \mapsto f$ intertwining.

The action of $G$ and the corresponding action of $\mathfrak{g}$ on $I'(q,r,s)$ have been calculated by Sepanski and Stanke \cite{Sepanski1}. We will record these results, since they are used in later sections.

\begin{proposition}\label{GroupActionStandard}
Let $f \in I'(q,r,s)$, $(g,\epsilon)\in G_2$ with $g=\left(\begin{smallmatrix}a & b \\ c & d\end{smallmatrix}\right)$, and $(v_1,v_2,w)\in H_{2n+1}$. Then,
\begin{subequations}
 \begin{align}
((g,\epsilon).f)(t,x)&=(a-ct)^{r-q/2} \epsilon(g^{-1}.(t+z)) e^{\frac{-s c \left\|x\right\|^2}{a-ct}} f\left(\frac{dt-b}{a-ct},\frac{x}{a-ct} \right) \label{SL2StdAction}\\
 ((v_1,v_2,w).f)(t,x)&=e^{-s(v_1\cdot v_2-2v_2\cdot x-t\left\|v_2\right\|^2+w)}f(t,x-v_1+tv_2). \label{H3StdAction}
\end{align}
\end{subequations}
Let $u\in O(n)$ then $(u.f)(t,x)=f(t,u^{-1}x)$.
\end{proposition}

\begin{corollary}\label{sl2ActionsStd}
The action of $\left(\begin{smallmatrix}a & b \\ c & -a\end{smallmatrix}\right)\in\mathfrak{sl}(2,\R)$ on $I'(q,r,s)$ is given by the differential operator
\begin{equation}\label{AlgActionNonComp}(ct-a)\sum_{j=1}^nx_j\partial_j+(ct^2-2at-b)\partial_t+(ra-cs\left\|x\right\|^2-rct). \end{equation}
An element $(u,v,w)\in \mathfrak{h}_{2n+1}$ acts on $I'(q,r,s)$ by the differential operator
$$-\sum_{j=1}^{n}u_j\partial_j+t\sum_{j=1}^{n}v_j\partial_j +s(w-2v\cdot x).$$
\begin{proof}It follows from differentiating the group actions on $I'(q,r,s)$. 
\end{proof}
\end{corollary}

\section{Casimir Operators}\label{Casimirs}

Let $$E_n=\sum_{j=1}^n x_j\partial_j$$ denote the Euler operator on $\R^n$, $$\Omega_{\mathfrak{sl}(2,\R)}=\frac{1}{2}h^2-h+2e^+e^-$$ denote the Casimir element in the universal enveloping algebra of $\mathfrak{sl}(2,\R)$, and $\Omega_{\mathfrak{so}(n)}$ denote the Casimir element for $\mathfrak{so}(n)$. In \cite{Sepanski1} it was shown that the element $\Omega$ in the universal enveloping algebra of $\mathfrak{g}$ defined by $$\Omega = 2\Omega_{\mathfrak{sl}(2,\R)}-\Omega_{\mathfrak{so}(n)}-r(r+2)$$ acts on $I'(q,r,s)$ as the differential operator $$\Omega'=-(2r+n)E_n+\left\|x\right\|^2(4s \partial_t + \Delta_n)$$ where $\Delta_n$ denotes the $n$-th dimensional Laplacian.

In particular, for $r=-n/2$, $\Omega$ acts on $I'(q,r,s)$ by $\Omega'=\left\|x\right\|^2(4s \partial_t + \Delta_n)$ so that,
$$\ker(\Omega'- 2\lambda)=\ker\Big(4s \partial_t + \Delta_n-\frac{2\lambda}{\left\|x\right\|^2}\Big).$$

\begin{remark}
In case $r=-n/2$ and $s=-1/2$ (respectively, $s=i/2$) the invariant subspace $\ker(\Omega'- 2\lambda)$ is contained in the space of solutions of the heat (respectively, Schr\"{o}dinger equation) with singular potential $\frac{2\lambda}{\left\|x\right\|^2}$. From here on, we let $r=-n/2$. 
\end{remark}

\section{The Compact Picture}\label{cpctPic}

We have given the explicit intertwining isomorphism between $I(q,-\frac{n}{2},s)$ and the non-compact picture $I'(q,-\frac{n}{2},s)$. In this section, we realize $I(q,-\frac{n}{2},s)$ in a way that will allow us to determine the $K_2\times O(n)$-weight vectors explicitly. To that end, we observe that the group $G_2$ has Iwasawa decomposition $G_2=K_2A\overline{N}$ and notice that the multiplication map $K_2 \times A \overline{N} \to G_2$  induces a diffeomorphism $G\cong (K_2\times X) \times ((A\overline{N}\times O(n))\ltimes W)$. Since $((A\overline{N}\times O(n))\ltimes W)\subset \overline P$, an element $\phi \in I(q,-\frac{n}{2},s)$ is completely determined by its restriction to $K_2\times X$.

Since $K_2 \cong S^1$ via a $4\pi$-periodic isomorphism, smooth functions on $K_2\times X$ can be realized as smooth functions on $S^1\times \R^n$.  In turn, these functions can be extended by periodicity to smooth functions on $\R^{n+1}$. Via the restriction map $\phi \to \phi |_{K_2\times X}$, $I(q,-\frac{n}{2},s)$  can be realized in $C^\infty(\R^{n+1})$ in the following way $\phi \mapsto F$ if and only if $$\phi([(g_\theta,\epsilon_\theta),(y,0,0)])=F(\theta,y).$$

We denote the image of this map by $I''(q,-\frac{n}{2},s)$. One can give this space a $G$-module structure that makes the map intertwining. Calculating the action of $(g_\theta,\epsilon_\theta)$ and considering the periodicity conditions, it is straightforward to show that 
\begin{equation}\label{defI''}I''(q,-\frac{n}{2},s)=\{F\in C^\infty(\R^{n+1}) | F(\theta+j\pi,(-1)^jy)=i^{-jq}F(\theta,y) \}. \end{equation}

Since we established an isomorphism between $I(q,-\frac{n}{2},s)$ and $I'(q,-\frac{n}{2},s)$, there exists an induced isomorphism from $I'(q,-\frac{n}{2},s)$ to $I''(q,-\frac{n}{2},s)$. This isomorphism is given by $F \mapsto f$ where 
\begin{equation}
f(t,x)=(1+t^2)^{-n/4}e^{\frac{st\left\|x\right\|^2}{1+t^2}}F(\arctan t ,x(1+t^2)^{-1/2}).
\end{equation} 
Equivalently, one can define $f \mapsto F$ by
\begin{equation}\label{CompactPicIsomorphism}F(\theta,y) =(\cos\theta)^{-n/2}e^{-s\left\|y\right\|^2\tan\theta}f(\tan\theta,y\sec\theta)
\end{equation}
for $\theta \in (-\frac{\pi}{2},\frac{\pi}{2})$. Here $F$ can be extended to $\theta \in \R$ by first using  continuity at $\theta = \pm \pi/2$ and then by using $F(\theta+j\pi,(-1)^jy)=i^{-jq}F(\theta,y)$.

Under this isomorphism, via the chain rule, we obtain
\begin{subequations}\label{DerivativesInCompPic}
\begin{gather}
\partial_t\leftrightarrow\frac{1}{2}(-y\sin2\theta \partial_y+\cos^2\theta\partial_\theta+2sy^2\cos2\theta-1/2r\sin2\theta) \label{PartialTCompact}\\
\partial_x\leftrightarrow 2sy\sin\theta+\cos\theta \partial_y. \label{PartialXCompact}
\end{gather}
\end{subequations}

Define a standard basis $\left\{\eta^+, \kappa, \eta^-\right\}$ of $\mathfrak{sl}_2(\C)$ given by $$\kappa = i(e^--e^+)$$ and $$\eta^\pm=1/2(h\pm i(e^++e^-)).$$ Applying equations \eqref{DerivativesInCompPic} to the action in Corollary \ref{sl2ActionsStd}, it can be shown that the $\mathfrak{sl}_2$-triple  $\left\{\eta^+, \kappa, \eta^-\right\}$ acts on $I''(q,-\frac{n}{2},s)$ by the differential operators
\begin{gather}
\kappa = i\partial_\theta \label{ActionKappa}\\
\eta^\pm  = \frac{1}{2}e^{\mp2i\theta}\left(-E_n \mp i\partial_\theta-(n/2\pm 2is\left\|y\right\|^2)\right). \label{ActionEtas}
\end{gather}
In the following corollary, we use these actions to calculate the action of $\Omega'$ on $I''(q,-\frac{n}{2},s)$.
\begin{proposition}
If $\Omega ''$ denotes the differential operator by which the central element $\Omega$ acts on $I''(q,-\frac{n}{2},s)$ then
$$\Omega''=\left\|y\right\|^2\left(4s\partial_\theta+4s^2\left\|y\right\|^2+ \Delta_n\right). $$
\end{proposition}

\section{$K_2\times O(n)$-types in $\ker(\Omega''-2\lambda)$}\label{Struc1}

Let $K=K_2\times O(n)$. The goal of this section is to write the $K$-types in $\ker(\Omega''-2\lambda)$ explicitly. Write $\mathcal H _k(\R^n)$ for the space of harmonic polynomials of homogeneous degree $k$ and $\mathcal H _k(S^{n-1})$ for the restriction of elements of $\mathcal H _k(\R^n)$ to $S^{n-1}$. The $O(n)$-finite vectors in $C^\infty(S^{n-1})$ are the harmonic polynomials on $S^{n-1}$. That is $$C^\infty(S^{n-1})_{O(n)-finite}=\bigoplus_k \mathcal H _k(S^{n-1}),$$ where $k\in \mathbb Z^{\geq 0}$ for $n\geq 3$, $k\in \mathbb Z$ for $n=2$, and $k\in \{0,1\}$ for $n=1$.  We decompose $0\neq y\in \R^n$ in polar coordinates as $y=\rho \xi$ with $\rho = \left\|y\right\|$ and $\xi \in S^{n-1}$. 

\begin{proposition}\label{GenFormKTypes}
The space of $K$-finite vectors in $I''(q,-\frac{n}{2},s)$ is the span of all functions of the form $$F(\theta,y)=e^{-im\theta/2}\psi(\rho)h_k(y),$$ (when $y\neq 0$) where $m\in\mathbb Z$, $\psi \in C^{\infty}(0,\infty)$, $h_k \in \mathcal H _k(\R^n)$, and $$m\equiv q+2k \mod 4$$
with $F(\theta,y)$ extending smoothly to $y=0$ and $\lim_{\rho\to 0}\rho^k \psi(\rho)$ bounded.
\begin{proof}
This result is proved by Sepanski and Stanke \cite{Sepanski1}.
\end{proof}
\end{proposition}

\begin{lemma}\label{CondOnKtypes}
If $m \in \mathbb Z$, $\psi \in C^2(\R)$, and $h_k \in \mathcal H_k(\R^n)$, then a function of the form $F(\theta, y)=e^{-im\theta/2}\psi(\rho)h_k(y)$ is in $\ker(\Omega''-2\lambda)$ if and only if $\psi$ is annihilated by the differential operator $$\mathcal D = \rho^2\partial_\rho^2+(n-1+2k)\rho \partial_\rho +4s^2 \rho^4-2ism\rho^2-2\lambda$$
\begin{proof}
The result follows from the following calculation: $$\rho^2\Delta_n(\psi h_k)= \rho^2\psi''h_k+\rho((n-1)h_k+2kh_k)\psi'$$ together with the fact that $$\Omega''-2\lambda=\rho^2\left(4s\partial_\theta+4s^2\rho^2+ \Delta_n-\frac{2\lambda}{\rho^2}\right).\qedhere$$
\end{proof}
\end{lemma}

For $n\geq 2$, let $$A_n:=\{l(n+2j)\ | \ l, j\in \mathbb{Z}\text { and } 1\leq l\leq j+1\}.$$ In preparation for writing the $K$-finite vectors explicitly and determining the conditions for their existence when $n\geq 2$, we state the following lemma.

\begin{lemma}\label{ParityCondLambda}
If $n\geq 2$ is even, then $A_n=2\mathbb{Z}\cap \mathbb Z^{\geq n}$. If $n\geq 3$ is odd, then $$A_n= \{\gamma= 2^{r}a \in \mathbb Z^{\geq 0}\ |\ (a,2)=1 \text{ and } a\geq n+2^{r+1}-2 \}$$
\begin{proof}
The case when $n$ is even is clear by letting $l=1$. Let $n\geq 3$ be odd and let $B_n$ be the right hand side of the equation. We want to show that $A_n = B_n$. Let $l=2^r b$ where $(b,2)=1$ and $b\geq 1$ then $l(n+2j)= 2^r b(n+2j)=2^ra$ where $a= b(n+2j)$. Now, we have $$b(n+2j)\geq b(n+2(l-1)) =b(n+2(2^rb-1))\geq n+2^{r+1}-2, $$ which shows $A_n\subset B_n$.

The other inclusion follows almost immediately by noticing $$2^r(1+2k)=2^r(n+(2k+1-n))\geq 2^r(n+2^{r+1}-2).$$ From this we see that $1\leq 2^r \leq \frac{2k-n+1}{2}+1$ and thus any elements in $B_n$ are in $A_n$.
\end{proof}
\end{lemma}

For convenience, we state the following definition. It will be used to describe the eigenvalues of $\Omega"$ for which the space of $K$-finite vectors is non-empty. We must remark that it should not be confused with the standard usage of the term admissible for a representation.
\begin{definition}
\begin{enumerate}
	\item For $n=1$ we say that  an \textbf{eigenvalue $\lambda$ is admissible} if and only if $\lambda$ is a triangular number. That is $\lambda =\frac{1}{2}l(l-1)$ for some $l\in \mathbb Z^{\geq 0}$. Denote this set by $A_1$.
	\item For $n \geq 2$ we say that an \textbf{eigenvalue $\lambda$ is admissible} if and only if $\lambda \in A_n$.
	\item Let $n\geq 3$ and $\lambda\in A_n$. We call the pair $(l,k)\in  \mathbb Z^{\geq 1}\times\mathbb Z^{\geq 0}$, \textbf{$\lambda$-admissible} if and only if $\lambda = l(2l+2k+n-2)$.
	\item Let $n=2$ and $\lambda$ be admissible. The pair $(l,k)\in \mathbb Z^{\geq 1}\times \mathbb Z$ is \textbf{$\lambda$-admissible} if and only if $\lambda = l(2l+2k+n-2)$.
	\item When $n=1$, a pair $(l,0)$ is $\lambda$-admissible iff $\lambda = \frac{l(l-1)}{2}$.
\end{enumerate}
\end{definition}
\begin{figure}
	\centering
		\includegraphics[width=4in]{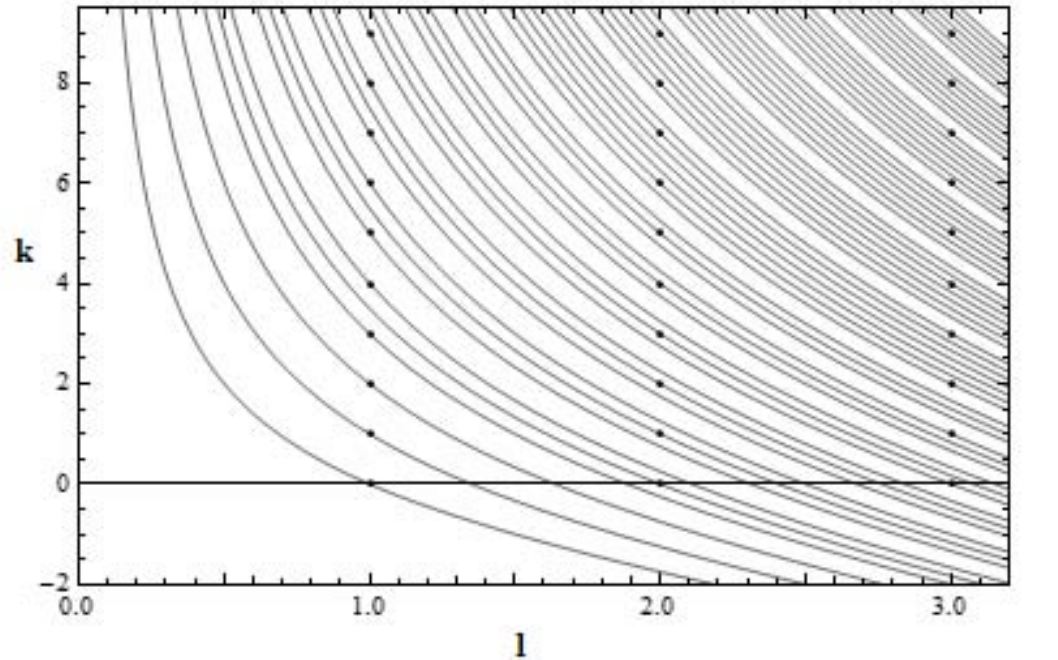}
	\caption{Admissible $\lambda$-level curves for $n=3$ in the $(l,k)$-plane.}
	\label{fig:lambdas}
\end{figure}

In Figure \ref{fig:lambdas}, the admissible $\lambda$ level curves for $n=3$ are represented. Notice that for $\lambda =0$, the corresponding level curve is the $k$-axis. The $\lambda$-admissible pairs $(l,k)$ are the points in $\mathbb Z^{\geq 1}\times \mathbb Z^{\geq 0}$ that intersect the $\lambda$ level curve. 

\begin{theorem}\label{KTypesForm}
\begin{enumerate}
 	\item \label{ItemTwo} The space of $K$-finite vectors in $\ker(\Omega''-2\lambda)\subset I''(q,-\frac{n}{2},s)$ is non-empty if and only if $\lambda$ is admissible.
 	\item \label{ItemThree} In this case, the space of $K$-finite vectors in $\ker(\Omega''-2\lambda)\subset I''(q,-\frac{n}{2},s)$ is spanned by the functions $F_{m,l,k}(\theta,y)$ of the form
 	\begin{equation*}\label{KFinForm}F_{m,l,k}(\theta,y)=e^{-im\theta/2}e^{-is\rho^2}\rho^{2l} h_k(y)_1F_1(\frac{m+4l+2k+n}{4},2l+k+n/2,2is\rho^2),\end{equation*} 
 	where the pair $(l,k)$ is $\lambda$-admissible and $m \equiv 2k+q \mod 4$.
\end{enumerate}
\begin{proof}
The one dimensional case was analyzed by Franco \cite{Franco} and the case when $\lambda=0$ has been studied by Sepanski and Stanke \cite{Sepanski1}. So, let $n\geq 2$ and $\lambda \neq 0$. The elements in $\ker(\Omega''-2\lambda)\subset I''(q,-\frac{n}{2},s)$ are of the form $F(\theta,y)=e^{-im\theta/2}\psi(\rho)h_k(y)$, they satisfy the conditions in Proposition \ref{GenFormKTypes}, and by Lemma \ref{CondOnKtypes}, $\mathcal D \psi =0$. Following Coddington \cite{Coddington}, the Frobenius method for this equation yields a solution space spanned by two linearly independent solutions. The indicial roots for this equation are $$r_\pm=\frac{1}{2}[2-2k-n\pm \sqrt{(n-2+2k)^2+8\lambda}].$$
Since $\mathcal D$ respects the decomposition of $\psi$ in terms of its even and odd components and these are determined by their value on $\R^{\geq 0}$ we may assume that $\psi$ is a function of $\rho^2$. Then, the first linearly independent solution is of the form $$\psi(\rho)=\rho^{r_+}(1+\sum_{j=1}^\infty c_j(r_+) \rho^{2j})$$ for some $c_j(r_+)\in \R$. This function extends to a smooth function of $y\in\R^n$ if and only if $r_+ \in 2 \mathbb Z^{\geq 0}$ or, equivalently, if and only if $$\lambda=l(2l+2k+n-2)$$ for $l\in  \mathbb Z^{\geq 0}$ where $r_+=2l$. The parity conditions of $\lambda$ stated in \eqref{ItemTwo} follow directly from this expression together with Lemma \ref{ParityCondLambda}. Moreover, it is clear that $l$ must be a divisor of $\lambda$. Solving for $k$ one obtains $k=\frac{\lambda}{2l}-l+1-\frac{n}{2}$. Then $k\in \mathbb Z$ iff $\lambda/l \equiv n \mod 2$ and, when $n\geq 3$, we need $k\geq 0$, which happens if and only if $1\leq l \leq 1/4(-(n-2)+\sqrt{(n-2)^2+8\lambda})$.

If $r_+ \in 2 \mathbb Z^{\geq 0}$, then $r_+-r_-\in Z^{\geq 0}$. If $r_+-r_-\neq 0$ then the second independent solution is of the form $\psi_2(\rho)=a\psi(\rho)\ln \rho + \rho^{r_-}(1+\sum_{j=0}^\infty c_j(r_-)\rho^{2j})$ and it is not continuous at zero because $r_-<0$. If $r_+-r_-=0$ then  $\lambda=0$ and $(k,n)=(0,2)$ is the only possible solution. In this case, it is known that $a=1$ so $\ln \rho$ makes the second solution not continuous at zero. Therefore, for each of these admissible pairs $(l,k)$, there exists a unique $K$-finite vector of the form $\psi(\rho)=\rho^{2l}(1+\sum_{j=1}^\infty c_j \rho^{2j})$.

To establish \eqref{ItemThree}, it now suffices to show that for fixed $(m,l,k)$, the corresponding $K$-finite vector is $\psi_{m,l,k}(\rho)=e^{-is\rho^2}\rho^{2l}\ _1F_1(\frac{m+4l+2k+n}{4},2l+k+n/2,2is\rho^2)$. Explicitly calculating $\mathcal D (\rho^{2l}e^{-is\rho^2}F(2is\rho^2))$, one obtains the following differential equation:
$$2is\rho^2F''(2is\rho^2)+(k+2l+\frac{n}{2}-2is\rho^2)F'(2is\rho^2)-
\frac{m+4l+2k+n}{4}F(2is\rho^2)=0.$$
Recall that the confluent hypergeometric differential equation is $$(z\partial_z^2+(b-z)\partial_z - a)_1F_1(a,b,z)=0$$ (Abramowitz and Stegun \cite{Abramowitz}). This equation has well known solutions in the form of confluent hypergeometric functions of the first and second kind. However, the smoothness condition required by being in $I''(q,r,s)$ shows that the unique solution corresponds to a multiple of the confluent hypergeometric function of the first kind. We may therefore take $F(2is\rho^2)=\phantom{.}_1F_1(\frac{m+4l+2k+n}{4},2l+k+n/2,2is\rho^2)$. As a function of $y$, this solution extends smoothly to a solution on $\R^n$.
\end{proof}
\end{theorem}

\section{Irreducible Subspaces of $\ker(\Omega''-2\lambda)$}\label{Struc2}
In this section, we look at the structure of $\ker(\Omega''-2\lambda)\subset I''(q,-\frac{n}{2},s)$ as an $\mathfrak{sl}_2\times O(n)$-module. To that end, we will explicitly compute the actions of the standard $\mathfrak{sl}_2$-basis. For these calculations, we will use the following properties of the confluent hypergeometric function (Abramowitz and Stegun \cite{Abramowitz}):
\begin{subequations}
\begin{gather}
\frac{d^n}{dz^n}\ _1F_1(a,b,z)=\frac{(a)_n}{(b)_n}\ _1F_1(a+n,b+n,z) \label{U0} \\
b\ _1F_1(a,b,z)-b\ _1F_1(a-1,b,z)-z\ _1F_1(a,b+1,z)=0 \label{U1} \\
\begin{split}
b ( 1-b+z)\,_1F_1(a,b,z)+b(b-1)\,_1F_1(a&-1,b-1,z) \\  & -az\,_1F_1(a+1,b+1,z)=0 \label{U2} \end{split} \\
\begin{split}
(a-1+z)\,_1F_1(a,b,z)+(b-a)\,_1F_1(a&-1,b,z) \\  & (1-b)\,_1F_1(a,b-1,z)=0 \label{U3} \end{split} \\
(a-b+1)\,_1F_1(a,b,z)-a\,_1F_1(a+1,b,z)+(b-1) \,_1F_1(a,b-1,z)=0 \label{U4}
\end{gather}
\end{subequations}
Combining \eqref{U1} with $a+1$ instead of $a$ and \eqref{U4} one obtains
\begin{equation}\label{Uno}
\,_1F_1(a,b,z)=\,_1F_1(a,b-1,z)-\frac{az}{b(b-1)}\,_1F_1(a+1,b+1,z).
\end{equation}
Using Equation \eqref{U4} with $b+1$ in place of $b$ and combining it with \eqref{U2}, one obtains
\begin{equation}\label{Dos}
\,_1F_1(a,b,z)=\,_1F_1(a-1,b-1,z)-\frac{b-a}{b-1}z\,_1F_1(a,b+1,z).
\end{equation}

\begin{theorem}\label{slActionKtypes}
For a $\lambda$-admissible pair $(l,k)$ with $h_k \in \mathcal H _k(\R^n)$ non-zero, for the $\mathfrak{sl}_2$ triple $\left\{\eta^+, \kappa, \eta^-\right\}$ we have
\begin{subequations}
\begin{gather}
\kappa.F_{m,l,k}(\theta,y)=\frac{m}{2}F_{m,l,k}(\theta,y), \\
\eta^\pm.F_{m,l,k}(\theta,y)=-\frac{\pm m+4l+2k+n}{4}F_{m\pm 4,l,k}(\theta,y).
\end{gather}
\end{subequations}
Lowest weight vectors occur if and only if $q\equiv n \mod 4$ and, in this case,
$$F_{(2k+4l+n),l,k}=e^{-\frac{1}{2}(2k+4l+n)i\theta}e^{is\rho^2}\rho^{2l}h_k$$
is a lowest weight vector. Highest weight vectors occur if and only if  $q+n\equiv 0\mod 4$ and, in this case,
$$F_{-(2k+4l+n),l,k}=e^{\frac{1}{2}(2k+4l+n)i\theta}e^{-is\rho^2}\rho^{2l}h_k$$
is a highest weight vector.

\begin{proof}
Let $p_m(\theta,\rho)=e^{-im\theta/2}e^{-is\rho^2}$. Since $\eta^\pm$ act by $-1/2e^{\mp2i\theta}[-E_n\mp i\partial_\theta+(n/2\mp 2is\rho^2)]$, we have \begin{multline*}\eta^+.F_{m,l,k}(\theta,y)=-p_{m+4}(\theta,\rho)\rho^{2l}h_k(y)\left(a_1F_1(a,b,z)+z\frac{a}{b}\phantom{x} _1F_1(a+1,b+1,z)\right)\end{multline*}
where $a=\frac{m+4l+2k+n}{4}$, $b=2l+k+n/2$, and $z=2is\rho^2$. Then, Equation \eqref{U1} implies
$$\eta^+.F_{m,l,k}(\theta,y)=-\frac{m+4l+2k+n}{4}F_{m+4,l,k}.$$
The action of $\eta^-$ is as follows:
\begin{multline*}\eta^-.F_{m,l,k}(\theta,y)=-1/2p_{m-4}(\theta,\rho)\rho^{2l}h_k(y)\left((2l+k-m/2+n/2-4is\rho^2)_1F_1(a,b,z)\right. \\ \left. +2z\frac{a}{b}\phantom{x} _1F_1(a+1,b+1,z)\right)\\ = -p_{m-4}(\theta,\rho)\rho^{2l}h_k(y)\left((b-a-z)_1F_1(a,b,z)\right. \\ \left. +z\frac{a}{b}\phantom{x} _1F_1(a+1,b+1,z)\right).\end{multline*}
Equation \eqref{U2} implies
$$\eta^-.F_{m,l,k}(\theta,y)=-p_{m-4}(\theta,\rho)\rho^{2l}h_k(y)\left((1-a)_1F_1(a,b,z) +(b-1) _1F_1(a-1,b-1,z)\right).$$
Equation \eqref{U4} implies
$$\eta^-.F_{m,l,k}(\theta,y)=-\frac{-m+4l+2k+n}{4}F_{m-4,l,k}.$$
The statement about the lowest and highest weights follows by observing that such vectors can occur only when $2k+n+4l \equiv m \mod 4$ or $2k+n+4l \equiv- m \mod 4$, respectively. This fact, together with the condition that $2k+q\equiv m\mod 4$, gives the desired result. The form of the highest weight vectors follows from directly calculating the weight vectors with weight $m=-(2k+n+4l)$. The form of the lowest weight vectors follows in the same way, but with $m=2k+n+4l$.
\end{proof}
\end{theorem}
\begin{definition}
Let $\ker(\Omega''-2\lambda)_K$ denote the $K$-finite vectors of $\ker(\Omega''-2\lambda)\subset I''(q,-\frac{n}{2},s)$. For a $\lambda$-admissible pair $(l,k)$ define $H_{l,k}\subset \ker(\Omega''-2\lambda)_K$ by $$H_{l,k}:=\spn\{F_{m,l,k} \ | \ m \equiv 2k+q \mod 4\}$$
If $q\equiv n \mod 4$, define $H_{l,k}^+\subset H_{l,k}$ by
$$H_{l,k}^+:=\{F_{m,l,k}\ | m\geq(2k+4l+n), \ m \equiv 2k+q \mod 4\ \}.$$
If $q\equiv -n \mod 4$, define $H_{l,k}^-\subset H_{l,k}$ by
$$H_{l,k}^-:=\{F_{m,l,k}\ | m\leq-(2k+4l+n), \ m \equiv 2k+q \mod 4\ \}.$$
\end{definition}

\begin{proposition}\label{irredLemma}
Let $\lambda$ be an admissible eigenvalue and $(l,k)$ a $\lambda$-admissible pair. Then, as $\mathfrak{sl}_2\times O(n)$-modules:
\begin{enumerate}
	\item If $q \not\equiv n\mod 4$ and $q \not\equiv -n\mod 4$, then $H_{l,k}$ is irreducible. Moreover, as an $\mathfrak{sl}_2\times O(n)$-module, $\ker(\Omega''-2\lambda)_K$ is decomposed as:
	$$\ker(\Omega''-2\lambda)_K=\bigoplus_{\substack{
\lambda\text{-admissible}\\
(l,k)}}H_{l,k}.$$
	\item If $q \equiv n\mod 4$ and $q \not\equiv -n\mod 4$, then $H_{l,k}^+$ is the unique irreducible $\mathfrak{sl}_2\times O(n)$-submodule of $H_{l,k}$.
	\item If $q \not\equiv n\mod 4$ and $q \equiv -n\mod 4$, then $H_{l,k}^-$ is the unique irreducible $\mathfrak{sl}_2\times O(n)$-submodule of $H_{l,k}$.
	\item If $q \equiv n\mod 4$ and $q \equiv -n\mod 4$, then $H_{l,k}^+$ and $H_{l,k}^-$ are the only irreducible $\mathfrak{sl}_2\times O(n)$-submodules of $H_{l,k}$. A composition series for $H_{l,k}$ is $$0\subset H_{l,k}^+ \subset H_{l,k}^+\oplus H_{l,k}^-\subset H_{l,k}.$$
\end{enumerate}
\begin{proof}
By Theorem \ref{slActionKtypes}, the representation is irreducible whenever $\pm(2k+4l+n)\neq m$ for any $m\equiv 2k+q \mod 4$, it has a highest or lowest weight submodule otherwise. 
\end{proof} 
\end{proposition}
\begin{remark}
The direct sum in Proposition \ref{irredLemma} is finite because, as a consequence of Proposition \ref{KTypesForm}, the set of $\lambda$-admissible pairs $(l,k)$ is finite for every admissible $\lambda$. 

When $H_{l,k}^+$ is non-empty, the representation is isomorphic to the $2k+4l+n$-th tensor product of the oscillator representation. Its dual occurs when $H_{l,k}^-$ is non-empty.
\end{remark}
\begin{figure}
	\centering
		\includegraphics{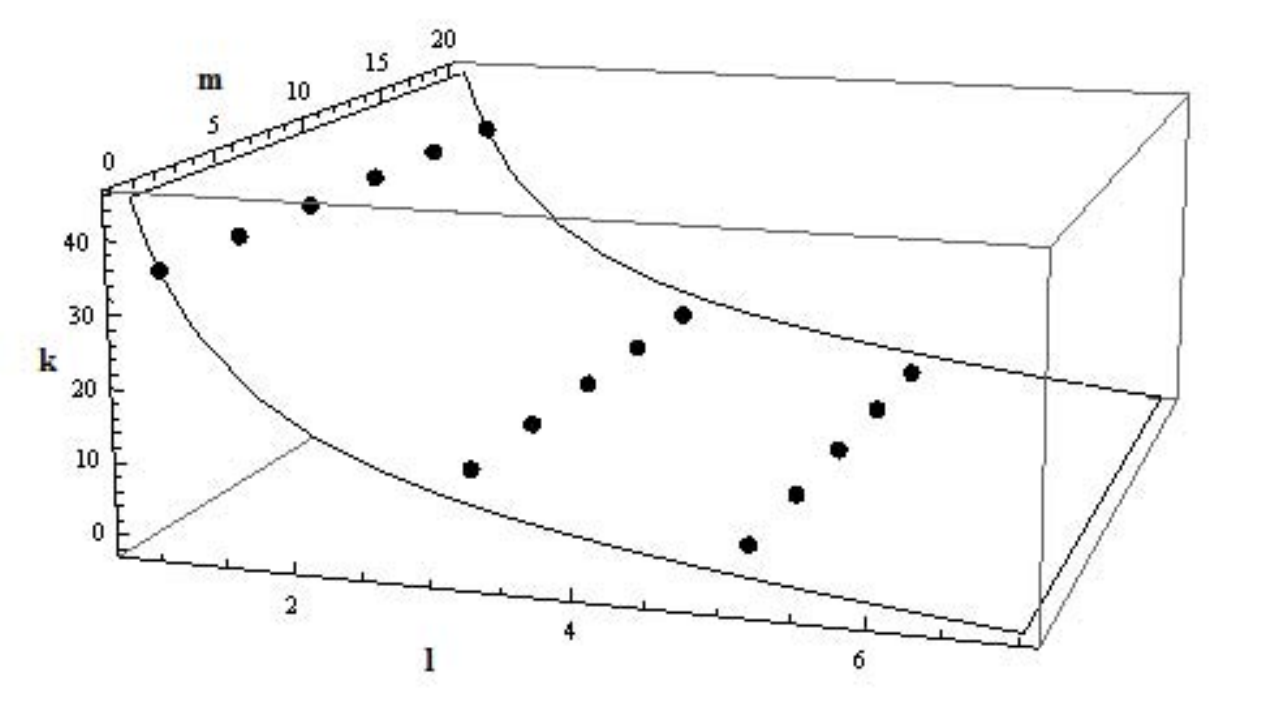}
	\caption{$K$-finite vectors of weight $\frac{m}{2}$ in $\ker(\Omega''-2\lambda)\subset I''(0,-3/2,s)$ for $0\leq m\leq 20$, $n=3$ and $\lambda=75$.}
	\label{fig:lambda2}
\end{figure}
Pictorially, a basis of $\ker(\Omega''-2\lambda)_K$ is depicted in Figure \ref{fig:lambda2}. There, we look at the case $n=3$ and $\lambda=75$. On the $(l,k)$-plane, the $\lambda$-level curve intersects the integral lattice on three different points: $(5,2)$, $(3,9)$, and $(1,36)$. These are all the $\lambda$-admissible pairs for this particular case. At each of these pairs, we have an $\mathfrak{sl}_2\times O(n)$-module represented by a line where each dot corresponds to the complex span of $F_{m,l,k}$. Notice that on the $m$ direction, the points are separated by jumps of $4$ units. This corresponds to the action of $\eta^{\pm}$.

\section{Heisenberg Action}
In this section, we will calculate the action of the Heisenberg algebra on the $K$-finite vectors. As it turns out, elements of the Heisenberg algebra $\mathfrak{h}_{2n+1}^\C$ will take a $K$-finite vector $F_{m,l,k}\in \ker(\Omega''-2\lambda)$ and map it to a linear combination of $K$-finite vectors associated to possibly other eigenvalues. The fact that a pair $(l,k)\in \mathbb{Z}^{\geq 1} \times \mathbb{Z}^{\geq 0}$ determines a unique $\lambda$ will be used to determine these eigenvalues explicitly. We begin with a lemma that was proved by Sepanski and Stanke \cite{Sepanski1} that will be used in the calculation of the actions below.

\begin{lemma}\label{HarmonicCondition}
Let $(k,n)\in \mathbb{Z}^{\geq 0}\times \mathbb{Z}^{\geq 0}$ and define constants $c_{k,n}=\frac{1}{2k+n-2}$ for $(k,n)\neq (0,2)$ and $c_{0,2}=0$. If $1\leq j\leq n$ and $h_k \in \mathcal H_k(\R^n)$, then
$$y_jh_k(y)-c_{k,n}\rho^2\in  \mathcal H_{k+1}(\R^n).$$
Moreover, if $h_k\neq 0$ and $(k,n)\neq (1,1)$, then there exists a $j\in \{1,2,...,n\}$ for which $y_jh_k(y)-c_{k,n}\rho^2(\partial_jh_k)(y)\neq 0$.
\end{lemma}

Let $\{e_j\}$ be the standard basis of $\C^n$ and define $E_j^\mp:=(\pm ie_j,e_j,0)\in \mathfrak{h}_{2n+1}^\C$. By Lemma \ref{sl2ActionsStd}, $E_j^\mp$ act on $I''(q,-\frac{n}{2},s)$ by $e^{\pm i\theta}(\mp i\partial_j-2sy_j)$.

\begin{proposition}\label{ActionEs}
For non-zero $F_{m,l,k}$,
\begin{multline}\label{Eplus}
E_j^+.F_{m,l,k}=2ilF_{m+2,l-1,k+1}-s\frac{(2l+2k+n-2)(m+2k+4l+n)}{2(k+2l+n/2-1)(k+2l+n/2)}F_{m+2,l,k+1}\\+\frac{i}{2k+n-2}\left((2l+2k+n-2)F_{m+2,l,k-1} \right.\\ \left.+\frac{2isl(m+2k+4l+n)}{2(k+2l+n/2-1)(k+2l+n/2)}F_{m+2,l+1,k-1}\right)\\ 
\end{multline}
and
\begin{multline}\label{Eminus}
E_j^-.F_{m,l,k}=-s\frac{(2-2l-2k-n)(4l+2k+n-m)}{2(2l+k+n/2-1)}F_{m-2,l,k+1}-2ilF_{m-2,l-1,k+1}\\-\frac{i}{2k+n-2}\left((2l+2k+n-2)F_{m-2,l,k-1}-isl\frac{4l+2k+n-m}{2(2l+k+n/2-1)}F_{m-2,l+1,k-1}\right).
\end{multline}

\begin{proof}
Using $a$ and $b$ as in the proof of Theorem \ref{slActionKtypes}, we explicitly calculate \begin{multline*}
E_j^+.F_{m,l,k}=ip_{m+2}(\theta,\rho)\left[ (2l\rho^{2(l-1)}h_ky_j+\rho^{2l}\partial_jh_k)_1F_1(a,b,z) \right. \\ \left. +4is\rho^{2l}y_jh_k\frac{a}{b}\phantom{x} _1F_1(a+1,b+1,z) \right].
\end{multline*}
Lemma \ref{HarmonicCondition} implies that there exists a possibly zero harmonic polynomial, $h_{k+1,j}\in \mathcal{H}_{k+1}(\R^n)$, such that $y_j h_k= h_{k+1,j}+c_{k,n}\rho^2\partial_jh_k$. Then we can write
\begin{multline*}
E_j^+.F_{m,l,k}=ip_{m-2}(\theta,\rho)\rho^{2(l-1)}\left(2l\,_1F_1(a,b,z)+4is\rho^2\frac{a}{b}\,_1F_1(a+1,b+1,z) \right)\\ \cdot h_{k+1,j}+ic_{k,n}\rho^{2l}p_{m-2}(\theta,\rho)\left((2l+c_{k,n}^{-1})\,_1F_1(a,b,z) \right. \\ \left. +4is\rho^2\frac{a}{b}\,_1F_1(a+1,b+1,z) \right)\partial_j h_{k}.
\end{multline*}
Using \eqref{Uno} we obtain, 
\begin{multline*}
E_j^+.F_{m,l,k}=i\rho^{2(l-1)}p_{m-2}(\theta,\rho)\left(2l\,_1F_1(a,b-1,z)\right. \\ \left. +\frac{az}{b}(2-\frac{l}{b-1})\,_1F_1(a+1,b+1,z) \right)h_{k+1,j}+ic_{k,n}\rho^{2l}p_{m-2}(\theta,\rho) \\ \cdot \left((2l+c_{k,n}^{-1})\,_1F_1(a,b-1,z) +\frac{az}{b}(2-\frac{l+c_{k,n}^{-1}}{b-1})\,_1F_1(a+1,b+1,z) \right)\partial_j h_{k},
\end{multline*}
which gives
\begin{multline*}
E_j^+.F_{m,l,k}=2ilF_{m+2,l-1,k+1}-s\frac{(2l+2k+n-2)(m+2k+2l+n)}{2(k+2l+n/2-1)(k+2l+n/2)}F_{m+2,l,k+1}\\+\frac{i}{2k+n-2}\left((2l+2k+n-2)F_{m+2,l,k-1}\right. \\ \left. +\frac{2isl(m+2k+2l+n)}{2(k+2l+n/2-1)(k+2l+n/2)}F_{m+2,l+1,k-1}\right). 
\end{multline*}
The action of $E^-$ is as follows:
\begin{multline*}
E_j^-.F_{m,l,k}=-i\rho^{2(l-1)}p_{m-2}(\theta,\rho)\left((2l-2z)y_jh_k+\rho^2\partial_jh_k \right)\,_1F_1(a,b,z)\\+2\rho^{2(l-1)}p_{m-2}(\theta,\rho)z\frac{a}{b}y_jh_k\,_1F_1(a+1,b+1,z).
\end{multline*}
Substitute $y_j h_k= h_{k+1,j}+c_{k,n}\rho^2\partial_jh_k$ and use Equation \eqref{U2} to obtain
\begin{multline*}
E_j^-.F_{m,l,k}=-i\rho^{2(l-1)}p_{m-2}(\theta,\rho)\left((2l+2-2b)\,_1F_1(a,b,z) \right. \\ \left. +2(b-1)\,_1F_1(a-1,b-1,z) \right)h_{k+1,j}-ic_{k,n}\rho^{2l}p_{m-2}(\theta,\rho) \\ \cdot \left((2l+2-2b+c_{k,n}^{-1})\,_1F_1(a,b,z)+2(b-1)\,_1F_1(a-1,b-1,z) \right)\partial_j h_{k}.
\end{multline*}
Using Equation \eqref{Dos} we obtain
\begin{multline*}
E_j^-.F_{m,l,k}=-i\rho^{2(l-1)}p_{m-2}(\theta,\rho)\left((2l+2-2b)\frac{b-a}{b-1}z\,_1F_1(a,b+1,z)\right. \\ \left.+l\,_1F_1(a-1,b-1,z) \right)h_{k+1,j}-ic_{k,n}\rho^{2l}p_{m-2}(\theta,\rho) \left((2l+2-2b+c_{k,n}^{-1})\right. \\ \left. \cdot\frac{b-a}{b-1}z\,_1F_1(a,b+1,z)+(2l+c_{k,n}^{-1})\,_1F_1(a-1,b-1,z) \right)\partial_j h_{k}.
\end{multline*}
Substituting for $a$ and $b$ yields the desired result.
\end{proof}
\end{proposition}

\begin{definition}\label{Def3}
Let $$H=\bigoplus_{\substack{
k\geq 0\\
l\geq 1}}H_{l,k}.$$
and $H_0=\bigoplus_{k\geq 0}H_{0,k}.$
If $n\equiv q \mod 4$, define
$$H^+=\bigoplus_{\substack{
k\geq 0\\
l\geq 1}}H^+_{l,k}$$
and $H^+_0=\bigoplus_{k\geq 0}H^+_{0,k}.$
If $n\equiv -q \mod 4$, define
$$H^-=\bigoplus_{\substack{
k\geq 0\\
l\geq 1}}H^-_{l,k}.$$
and $H^-_0=\bigoplus_{k\geq 0}H^-_{0,k}.$
\end{definition}

\begin{figure}
	\centering
		\includegraphics{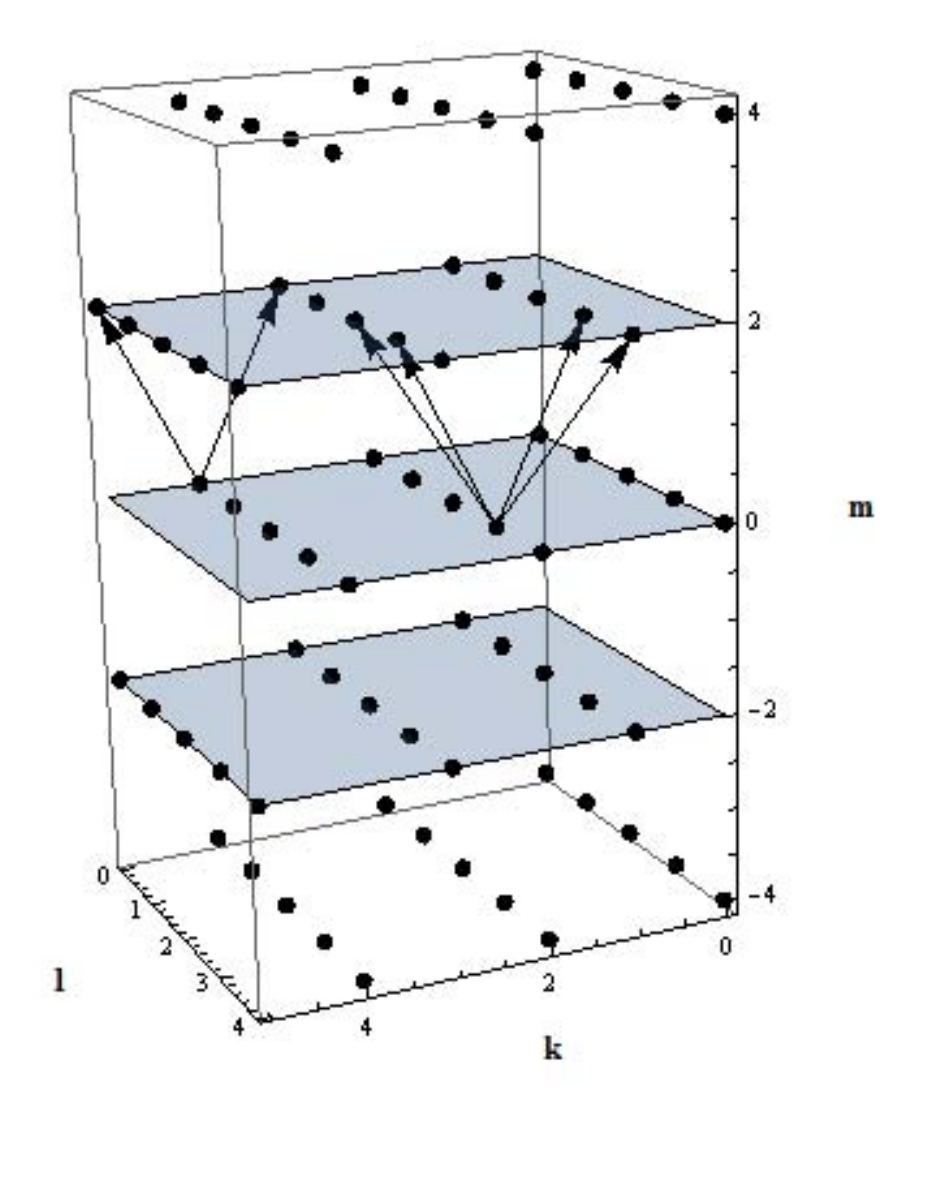}
	\caption{Action of $E_j^+$ on $H_0\oplus H$}
	\label{fig:Eplus}
\end{figure}

\begin{remark}
In Definition \ref{Def3}, the sums are defined with $k\geq 0$. This is the case when $n\geq 3$. However, when $n=2$ the sum must be taken over $k\in\mathbb Z$.
\end{remark}

In Figure \ref{fig:Eplus} we show how the element $E_j^+$ in the Heisenberg algebra acts on the space of $K$-finite vectors.

\begin{theorem}\label{BigTheorem}
As $\mathfrak g$-modules:
\begin{enumerate}
\item\label{BR1} If $q \not\equiv n\mod 4$ and $q \not\equiv -n\mod 4$ then $H_0$ is the unique irreducible submodule of $H$.
	\item\label{BR2} If $q \equiv n\mod 4$ and $q \not\equiv -n\mod 4$, then a composition series for $H$ is
	$$0\subset H_0^+\subset H_0 \subset H_0 \oplus H^+ \subset H$$
	\item \label{BR3}If $q \not\equiv n\mod 4$ and $q \equiv -n\mod 4$, then a composition series for $H$ is
	$$0\subset H_0^-\subset H_0 \subset H_0 \oplus H^- \subset H$$
	\item\label{BR4} If $q \equiv n\mod 4$ and $q \equiv -n\mod 4$, then a composition series for $H$ is
	$$0\subset H_0^-\subset H_0^+\oplus H_0^-\subset H_0 \subset H_0 \oplus H^- \subset H_0 \oplus H^- \oplus H^+ \subset H$$
\end{enumerate}
\begin{proof}
The statement in item \eqref{BR1} follows by noticing that, by Proposition \ref{ActionEs}, when $l=0$, the terms where the parameter $l$ is changed are annihilated by the Heisenberg algebra. This, together with the fact that $H_0$ is irreducible under the action of $\mathfrak g$ (Sepanski and Stanke \cite{Sepanski1}), gives the result.

The proofs of \eqref{BR2} and \eqref{BR3} are essentially the same. Therefore, we only  look at \eqref{BR2}. The first two inclusions in the composition series are a consequence of Proposition \ref{irredLemma} when $l=0$. In order to show the irreducibility of $H_0 \oplus H^+/H_0$ one has to notice that the actions of $E^\pm_j$ ``respects" the highest weight structures. More precisely,  Proposition \ref{ActionEs} implies that $$E_j^-.F_{(2k+4l+n),l,k}=-2ilF_{(2k+4l+n)-2,l-1,k+1}-i\frac{2l+2k+n-2}{2k+n-2}F_{(2k+4l+n)-2,l,k-1}$$
and this is a linear combination of highest weight vectors. In the same way, it can be seen from \eqref{Eplus} that $E_j^+$ maps a highest weight vector to a linear combination of elements in the highest weight modules corresponding to the triples $(m+2,l-1,k+1)$, $(m+2,l,k-1)$, $(m+2,l,k+1)$, and $(m+2,l+1,k-1)$. The elements corresponding to the first two triples are highest weight vectors and the latter get mapped to one by the action of $\eta^+$. The rest of the composition series in \eqref{BR2} is clear. 
\end{proof}
\end{theorem}

\begin{remark}
Suppose that $(l,k)$ is a $\lambda$-admissible pair. Then, the action of $E_j^\pm$ sends $F_{m,l,k}$ to a linear combination of $F_{m\pm2,l-1,k+1}$, $F_{m\pm2,l,k-1}$, $F_{m\pm 2,l,k+1}$, and $F_{m\pm 2,l+1,k-1}$. However the pairs $(l+1,k-1)$,  $(l-1,k+1)$, $(l,k-1)$, and $(l,k+1)$ are, in general, not admissible for $\lambda$, but they are admissible for different eigenvalues.

Therefore, $K$-finite vectors in $\ker(\Omega''-2\lambda)$ get sent, by $E_j^\pm$, to a linear combination of $K$-finite vectors in $\ker(\Omega''-2(\lambda\pm(2k+2l+n-2)))$ and in $\ker(\Omega''-2(\lambda\pm l))$.
\end{remark}


\end{document}